\documentclass[11pt]{article}
\usepackage{latexsym}
\usepackage{epsfig}
\usepackage{amsfonts}
\usepackage{amsmath}
\usepackage{amssymb}
\usepackage[latin1]{inputenc}

\newtheorem{definition}{Definition}
\renewcommand{\leq}{\leqslant}
\newcommand{\soussection}[1]{\vskip5mm\noindent\textbf{#1}\par}

\begin{document}

\title{An optimal algorithm to generate tilings}

\author{S\'ebastien Desreux 
	\thanks{
		LIAFA, UMR 7089 CNRS, Universit\'e Paris 7 Denis Diderot, 
		Case 7014, 2, place Jussieu, F-75251 Paris cedex 05, France, 
		\texttt{Sebastien.Desreux@liafa.jussieu.fr}
	} 
	\and Eric R\'emila 
	\thanks{
		LIP, umr 5668 CNRS-INRIA,  ENS Lyon, 46 all\'ee d'Italie, 
		69364 Lyon cedex 07, 
		France
		and GRIMA,
		IUT Roanne, 20 avenue de Paris, 42334 Roanne cedex, 
		France \texttt{Eric.Remila@ens-lyon.fr}
	}
}

\date{}
\maketitle

\begin{abstract}

A lot of progress has been made in tiling theory in the last ten years after
Thurston (\cite{Thu90}), building on previous work by Conway and Lagarias
(\cite{CL90}), introduced height functions as a tool to encode and study
tilings.

This allowed the authors of this paper, in previous work (\cite{Rem99},
\cite{Des01}), to prove that the set of lozenge (or domino) tilings  of
a hole-free, general-shape domain in the plane can be endowed with a
distributive lattice structure. 

In this paper, we see that this structure allows us in turn to construct an
algorithm that is optimal with respect to both space and execution time to
generate all the tilings of a domain~$D$. We first recall some results about
tilings and then we describe the algorithm.

\end{abstract}

\section{Background}

\soussection{Tiles}

The plane is endowed with either the square or triangular regular lattice
whose cells are colored black and white as on a chessboard. This induces a
direction on the edges of the lattice: They are directed clockwise around
black cells and (consequently) counterclockwise around white cells.  A {\em
domain} or \emph{region} is a finite and simply connected union of cells of
the lattice. The boundary of a domain $D$ will be denoted by $\partial D$.

A {\em lozenge} (resp. \emph{domino}) is a union of two cells of the
triangular (resp. square) grid sharing an edge, which is called the {\em
central axis} of the lozenge (resp. domino). This yields two possible shapes
for dominoes (vertical or horizontal) and three shapes for lozenges. A
\emph{domino tile} (resp. \emph{lozenge tile}) is a domino (resp. lozenge) of
either shape. A {\em tiling} of a domain is a set of tiles that cover the
whole area with neither gap nor overlap.

\soussection{Height functions} 

The {\em height functions}, introduced by W. P. Thurston (\cite{Thu90}) and
independently in the statistical physics litterature (see \cite{Bur} for a
review) and precisely studied and generalized by several authors (\cite{Cha},
\cite{Propp}, \cite{Rem99}, \cite{Des01}) are a very powerful tool to study
tilings.  A lozenge tiling $T$ of a domain $D$ can be encoded by a height
function $h_{T}$ defined as follows: Fix an origin vertex $O$ on the boundary
of $D$ and set $h_T(O) = 0$; if $(v, v')$ is a directed edge such that
$[v\,;\, v']$ is the central axis of a lozenge of $T$, then $h_T(v') = h_T(v)
- 2 $; otherwise, $h_T(v') = h_T(v) + 1 $. This definition is coherent since
it is coherent for each triangle and $D$ is simply connected. Similarly, for
dominoes, $h_T(v') = h_T(v) - 3 $ if $[v\,;\,v']$ is the central axis of a
domino, $h_T(v') = h_T(v) + 1 $ otherwise.

Height functions encode tilings: Only one such function is associated to a
tiling and a tiling can be reconstructed from a height function by drawing
only the edges whose endpoints have a height difference of~1.

\soussection{Lattice structure}

Let $(T, T')$ be a pair of tilings of $D$. We say that $T\leq T'$ if $h_T(v)
\leq h_{T'}(v)$ for each vertex $v$ of $D$.  The functions $h_{\text{inf}(T,
T') }= \text{min}(h_T, h_{T'})$ and $h_{\text{sup}(T, T')} = \text{max}(h_T,
h_{T'})$ are themselves height functions that encode tilings (\cite{Rem99},
\cite{Des01}), which implies that the set of the tilings of $D$ has a
structure of distributive lattice (see for instance \cite{Dav-Pri} for an
introduction to lattice theory).

\soussection{Flips} 

\begin{figure}
\begin{center}
	\epsfig{file=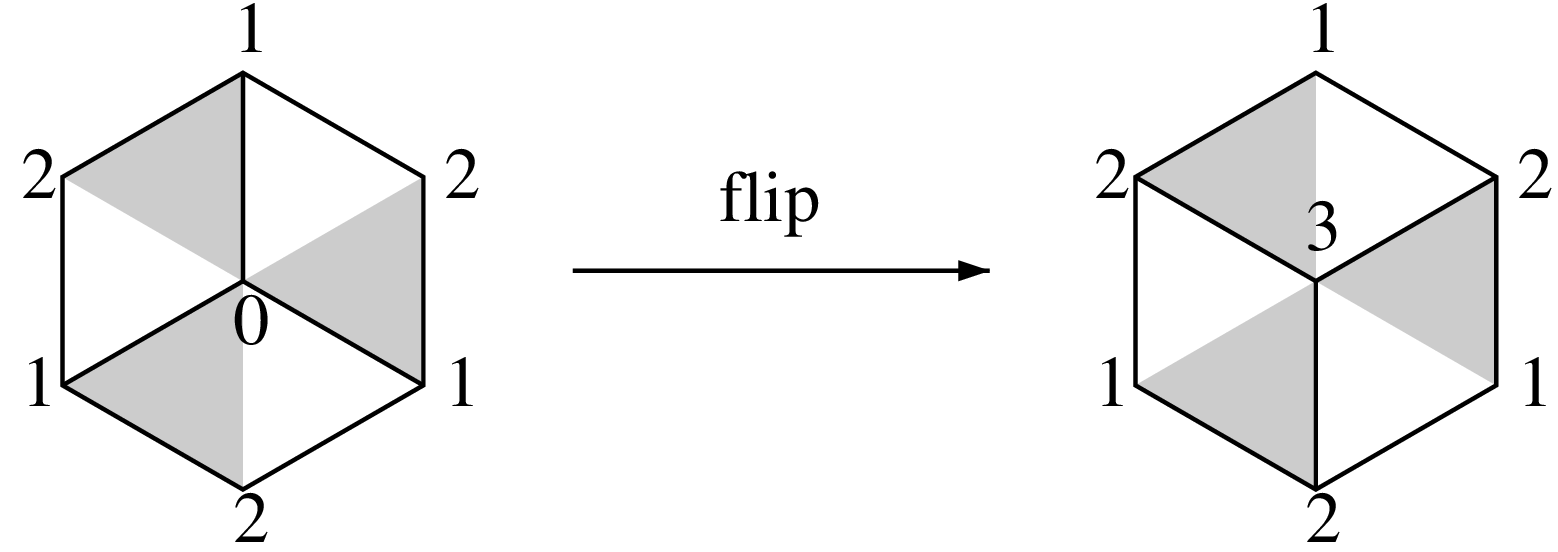, width=7cm}
	\vskip5mm
	\mbox{}\hskip2mm\epsfig{file=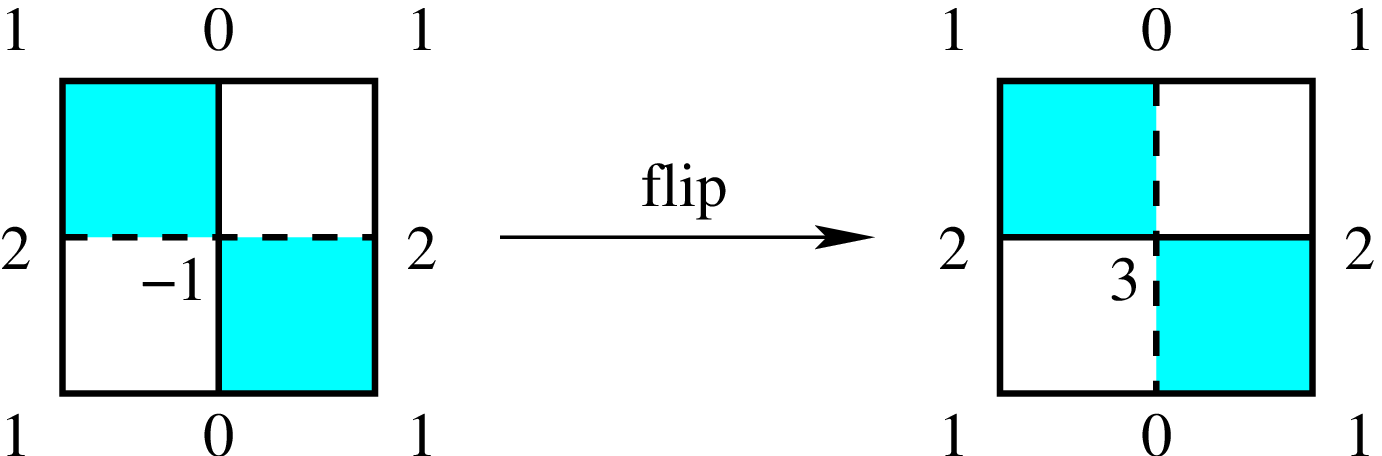, width=7cm}
	\label{flips}
	\caption{Flips in the triangular et square grids}
\end{center}
\end{figure}

Let $v$ be a vertex in the interior of $D$ such that all the directed edges
ending in $v$ are central axes of lozenges (resp. dominoes) in a tiling $T$. A
flip is the replacement of these three lozenges (resp. two dominoes) by three
lozenges (resp. two dominoes) whose central axis are edges starting in $v$.  A
flip transforms a local minimum of the height function into a local maximum.
See Figure~\ref{flips}.

A new tiling $T_{\text{flip}}$ is thus obtained; $T$ and $T_{\text{flip}}$ are
comparable for the order defined above. More generally, $T \leq T'$ if and
only if there exists an increasing sequence $(T = T_{0}, T_{1}, \ldots, T_{p}
= T')$ of tilings such that $\forall\ 0 \leq i < p$, $T_{i+1}$ is deduced from
$T_{i}$ by a flip. As a corollary we have the {\em flip connectivity}: Given
any pair $(T, T')$ of tilings of $D$, one can pass from $T$ to $T'$ by a
sequence of flips and, more precisely, the minimal number of flips to pass
from $T$ to $T'$ is $\sum_{v} |h_T(v) - h_{T'}(v)| / \lambda$, where
$\lambda=3$ for lozenges and $\lambda=4$ for dominoes.

\soussection{Thurston's algorithm}

There exists a minimal tiling whose corresponding height function has no local
maximum except on the boundary of $D$; indeed, a downward flip could otherwise
be performed on the local maximum, yielding a new minimal tiling. From this
property one deduces a linear algorithm which constructs the minimal tiling if
$D$ can be tiled, or proves that $D$ is not tileable (\cite{Thu90}).

First, a vertex on $\partial D$ must be selected and given an arbitrary
height, usually~0. Then the heights of all the vertices on the boundary of $D$
follow and they have the same height in all the tilings of $D$. Since the
local maxima of the height function must lie on $\partial D$, let us
select one such vertex and place tiles that cover it. There is
only one way to proceed without introducing a local maximum in the interior of
$D$. One can apply the same procedure to the remaining domain and a tiling is
thus built if at all possible.

This tiling is the minimal element of the lattice of the tilings of $D$. A
symmetric construction yields the maximal element, in which the height
function has no local minimum except on the boundary of $D$.

\section{Generalized Thurston algorithm}

Thurston's algorithm allows one to construct a particular tiling of a domain
$D$, using the fact that there exists a tiling whose height function has no
local maximum except on $\partial D$. But ordinary tilings do have height
functions which admit local maxima in the interior of $D$: Is there a way to
generalize Thurston's algorithm so that it can construct any tiling of $D$?

We will first reinterpret Thurston's algorithm with Birkhoff's representation
theorem for finite distributive lattices, then exhibit a generalized version
and finally give an example.

\soussection{Link with Birkhoff's representation theorem}

Let $T$ denote any tiling of a fixed domain $D$. Let $V$ denote the set of the
vertices in the interior of $D$ and $S \subset V$ the vertices on which the
height function $h_T$ reaches a local maximum. A downward flip can be applied
on any element of $S$ but on no element of $V \setminus S$. Therefore, $T$ is
a minimal tiling with respect to the heights on the vertices in $S$. In order
to characterize $T$, it is enough to know $h_T(v)$ for every $v \in S$. In
lattice terms, $T$ is the infimum of all the tilings of $D$ which have fixed
values on the elements of $S$.

What does the set $S$ represent? First, let us suppose that $S$ contains only
one vertex $v$. Then $T$ can be obtained by an upward flip from only one other
tiling, namely the one obtained by applying a downward flip on $v$ in $T$.
This property defines a meet-irreducible element of the lattice of the tilings
of $D$. Is $S$ contains more than one vertex, it can be viewed as a collection
of meet-irreducible elements. Birkhoff's representation theorem (see for
instance \cite{Dav-Pri}) allows us to formalize this idea:

\bigskip \noindent\textbf{\textit{Birkhoff's representation theorem}}

\noindent \textit{Any finite distributive lattice is isomorphic to the lattice
of the ideals of the order of its meet-irreducible elements.}

\bigskip

In other words, it is legitimate to regard a tiling $T$ as a collection (down
set) of meet-irreducible elements. These admit a simple characterization in
the case of tilings: A tiling is a meet-irreducible element of the lattice if
and only if its height function admits exactly one local maximum in the
interior of $D$.

From this point of view, the minimal element of the lattice, which is
precisely the tiling constructed by Thurston's algorithm, corresponds to the
empty ideal.

\soussection{Generalized algorithm}

We now undertake to construct any tiling of $D$ using Thurston's idea to cover
vertices so that no local maximum of the height function can appear.

\medskip

A meet-irreducible element of the lattice of the tilings of $D$ is
characterized by a unique pair $(v,h(v))$ where $v$ is the only vertex in the
interior of $D$ on which the height function admits a local maximum and $h(v)$
is the value of this local maximum.

Let $S$ be a set of vertices in the interior of $D$ and let $\mathcal{S}$ be
constituted of the pairs $(v,h(v))$ where $v \in S$ and $h(v)$ is any number.
According to Birkhoff's theorem, any tiling of $D$ is characterized by a set
$\mathcal{S}$. Our construction also allows sets $\mathcal{S}$ which
correspond to no tiling because $h$ may vary too rapidly along an edge.
We now suppose that $\mathcal{S}$ is a fixed set such that there exists at
least one tiling characterized by it.

If there exists at least one tiling whose height function coincides with
$h(v)$ for all $v \in S$, then there exists a smallest such tiling: It is the
infimum of the tilings satisfying this property. In particular, the height 
function of this smallest tiling, let us note it $T_{\mathcal{S}}$\,, can't
have a local maximum outside $S$, otherwise this would contradict the
minimality of $T_{\mathcal{S}}$\,. We can therefore apply Thurston's idea in 
order to construct this tiling:

\bigskip

\noindent\textbf{\textit{Generalized Thurston algorithm}}

\begin{itemize}

\item \textit{\textbf{Input}: A domain $D$, a vertex $v$ of $\partial D$ whose
height is~0, a subset $S$ of the vertices in the interior of $D$ and for each
$v \in S$, an integer $h(v)$.}

\item \textit{\textbf{Initialization}: Compute the height function of
$\partial D$. If a vertex receives two distinct heights, then $D$ is not
tileable.}

\item \textit{\textbf{Repeat}: Let $v \in \partial D \cup S$ be a vertex on
which the height function admits a global maximum. Place a tile whose frontier
covers $v$ in the only way that does not create a local maximum of the height
function. Remove $v$ from $D$ and update $\partial D$.}

\item \textit{\textbf{Until}: $D$ is tiled or one of its vertices was given
two different heights, in which case $D$ is not tileable.}

\end{itemize}

As the original, this algorithm uses each cell of the domain only once so it
is linear in the size of $D$. The space required is the one needed to store a
height function, which is $|D| \ln |D|$, where $|D|$ denotes the number of
cells of $D$.

\soussection{Example}

\begin{figure}[ht]
\begin{center}
        \epsfig{file=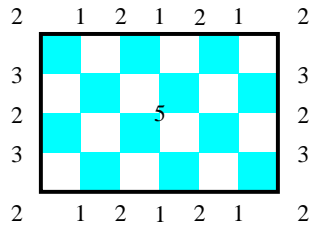, width=4cm}
        \caption{Generalized Thurston algorithm: step 1}
        \label{fig:TG:1}
\end{center}
\end{figure}

Let us illustrate the former algorithm with an example in the case of
dominoes. Our domain $D$ is a $6×4$ chessboard (see Figure~\ref{fig:TG:1});
the vertex at the center of $D$ has height~5. As in the case of Thurston's
original algorithm, one can fix arbitrarily the height of one vertex on
$\partial D$ and then compute the heights of all the vertices on $\partial D$.
Moreover, it is equivalent to proceed one tile at a time or to proceed in one
step for all the tiles covering a vertex of maximal height.

Initially, the height function admits a global maximum at the center of the
domain. This vertex must be covered by two dominoes, which can be either
vertical or horizontal. The latter case would yield a vertex having a
height greater than 5 so the former is the only possibility, otherwise we
would construct a tiling satisfying the conditions but that is not minimal in
this respect. We thus add two vertical dominoes and we update the height
function on their boundaries.

\begin{figure}[ht]
\begin{center}
        \epsfig{file=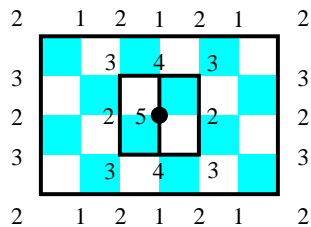, width=4cm}
        \caption{Generalized Thurston algorithm: step 2}
        \label{fig:TG:2}
\end{center}
\end{figure}

\begin{figure}[ht]
\begin{center}
        \epsfig{file=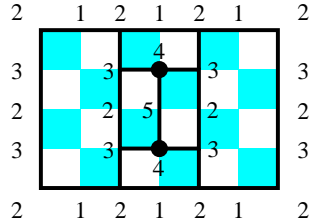, width=4cm}
        \caption{Generalized Thurston algorithm: step 3}
        \label{fig:TG:3}
\end{center}
\end{figure}

The maximal value of the height function is now~4 (see Figure~\ref{fig:TG:2}).
The geometry of the domain compels us to add horizontal dominoes, but let us
forget geometry and trust the algorithm. We can either add two vertical
dominoes or one horizontal one. The former case would attempt to add
a new local maximum, so only the latter is admissible. We thus add two
horizontal dominoes, as in Figure~\ref{fig:TG:3}.

\begin{figure}[ht]
\begin{center}
        \epsfig{file=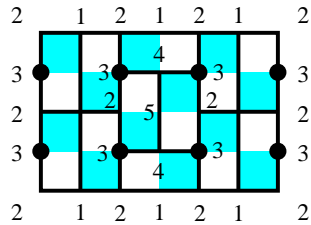, width=4cm}
        \caption{Generalized Thurston algorithm: step 4}
        \label{fig:TG:4}
\end{center}
\end{figure}

In the last step, we proceed simultaneously all the vertices of height~3, with
the same reasoning as above, and we obtain a complete tiling of $D$ (see
Figure~\ref{fig:TG:4}).

We see that the height function of the tiling thus constructed has no local
maximum except on the boundary and the center of $D$; it is thus the smallest
tiling satisfying the initial requirements.

If the conditions given initially had not corresponded to a tiling, a vertex
would have been given two different heights.

\section{Encoding of a tiling by a word}

In order to generate all the tilings of a domain, a natural idea is to encode
tilings by words and then to provide a way to find the successor of an element
in the lexicographic order. In this section, we examine a way to encode
tilings by words; the next section will provide a successor function.

\medskip

Since height functions encode tilings, it suffices to fix an arbitrary order
on the vertices of~$D$ (any order will do) to obtain an encoding of the
tilings by words. This encoding follows closely the height functions and is
therefore very natural. 

We will use a slightly different construction, although the one mentionned
above does work: Instead of the height function, we use a normalized height
function, defined as follows:

\[
	H(v) = \dfrac{h(v) - h_{\text{min}}(v)}{\lambda}
\]

\noindent where $h(v)$ is the value of the height function on the vertex $v$,
$h_{\text{min}}$ is the height function of the minimal tiling and $\lambda$ is
a normalization parameter, equal to~3 in the case of lozenges and to~4 in the
case of dominoes. This simple homography has two advantages: First, it unifies
the description of tilings by lozenges and by dominoes; second, the normalized
height is closely connected to flips. Indeed, $H(v)$ is the number of times a
flip has been applied to $v$ in any upward path going from the minimal tiling
to the tiling considered. And, of course, there is a one-to-one correspondence
between height functions and normalized height functions.

Moreover, it is convenient to represent a tiling by its phase:

\begin{definition}[Phase space]
Let $D$ be a tileable domain, $V$ the set of the vertices in the interior of
$D$ and $\varphi$ a numerotation function from $V$ to $\{ 1,\ldots, |V|\}$.

The \emph{phase} of a tiling $T$ is the set of pairs $\big(\, \varphi(v)\,,\,
H(\varphi(v))\,\big)$ where $v \in V$:

\[
	\Phi(T) = \big\{\ 
		\big(\, \varphi(v)\,,\,H(\varphi(v))\,\big)
		\ |\ v \in V
		\ \big\}
\]

The \emph{phase space} associated with a domain $D$ is the union of the
$\Phi(T)$ when $T$ runs through the set $\mathcal{T}(D)$ of the tilings of
$D$:

\[
	\Psi(D) = \bigcup\limits_{\mathcal{T}(D)} \Phi(T)
\]

\end{definition}

\begin{figure}[ht]
\begin{center}
        \epsfig{file=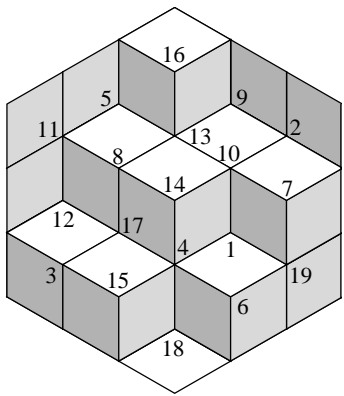, width=4.5cm}  
        \caption{An arbitrary order on the inner vertices of a domain}
        \label{fig:numerotation}
\end{center}
\end{figure}

We now give an example. Let us consider the domain of
Figure~\ref{fig:numerotation}, on which we have added an arbitrary order on
the inner vertices.

The normalized heights can be easily computed, but in the case of lozenges
they can also be read directly on the drawing. Indeed, since the lozenge group
is isomorphic to $\mathbb{Z}³$ (\cite{Thu90}), the result of a flip is,
visually, to add a cube. Take for instance vertex~16: one cube has been added
when starting from the minimal tiling of the domain, so its normalized height
is~1.

\begin{figure}[ht]
\begin{center}
        \epsfig{file=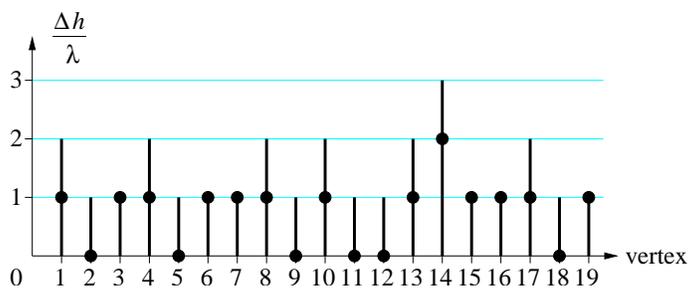, width=9cm}  
        \caption{The phase of a tiling}
        \label{fig:representation}
\end{center}
\end{figure}

In Figure~\ref{fig:representation}, we have represented the normalized heights
against the vertices' numbers. The dots correspond to the phase of the
tiling of Figure~\ref{fig:numerotation}. The vertical segments correspond to
the values that the normalized heights of a vertex can take in a tiling. Its
maximal value can be computed with Thurston's original algorithm in the
maximal tiling version.

\medskip

The encoding of the tiling follows easily from the phase diagram: it suffices
to read each normalized height in the order of the numerotation function:

\[
	µ(T) = H(1) \ \cdot\ H(2) \ \cdots\ H(|V|)
\]

For the example of Figure~\ref{fig:numerotation}, this yields the following
word:

\begin{figure}[ht]
\begin{center}
        \epsfig{file=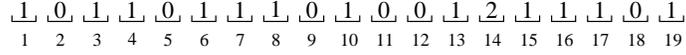, width=9cm}  
        \caption{Encoding of the tiling of Figure~\ref{fig:numerotation}}
        \label{fig:codage:0}
\end{center}
\end{figure}

\section{Exhaustive generation}

We have seen how to encode a tiling by a word. Since the lexicographic order
is a linear extension of the lattice structure, generating all the tilings of
a domain amounts to exhibiting a successor function.

\soussection{Successor of a word}

Let $w$ be an encoding of a tiling $T$ of $D$. The successor of $w$ in the
lexicographic order need not encode a tiling itself. We call the
\emph{successor} of $w$, and we will denote it by $\text{s}(w)$, the smallest
of the words greater than $w$ and that encode a tiling. We will show how to
construct it in two steps.

\medskip

Since tilings are connected by flips, there exists a series of tilings
going from $w$ to $\text{s}(w)$. Since $\text{s}(w)$ is minimal, the series
must contain exactly one upward flip. The first step is thus to determine the
local minima of the height function associated with $w$ and to select the
right-most position of the corresponding vertices in $w$. Applying one upward
flip yields a word $w'$ which differs from $w$ on only one position, let us
call it $i_0$:

\[
	\begin{cases}
	w'[i_0] = w[i_0] + 1 \\
	w'[i] = w[i] & \text{for $i \neq i_0$} \\
	\end{cases}
\]

$w'$ is greater than $w$ but it may be greater than $\text{s}(w)$. Consider
for instance the example of Figure~\ref{fig:numerotation} again. An upward
flip can be performed on vertex~18, which yields $\text{s}(w)$. At the next
iteration, however, the right-most candidate is~12. If an upward flip is
performed on this vertex, one obtains a tiling that is greater than
$\text{s}(\text{s}(w))$ since a downward flip could be performed on~18.

\medskip

The second step in finding the successor lies in the use of the generalized
Thurston algorithm. Suppose an upward flip has been performed (starting from
$w$) on the vertex $i_0$\,, yielding a word $w'$. The successor of $w$ has the
same values as $w'$ on positions 1 to $i_0$ included, but possibly smaller
ones for positions $> i_0$\,. It is indeed the smallest of the word that
coincide with $w'$ on the positions 1 to $i_0$\,.

We use the generalized algorithm by feeding it the heights already computed
for the vertices 1 to $i_0$ and letting it compute all the remaining ones. At
least one tiling exists under these conditions (the one associated with $w'$)
so the algorithm effectively yields a tiling, which has all the desired
properties.

\medskip

\begin{figure}[ht]
\begin{center}
        \epsfig{file=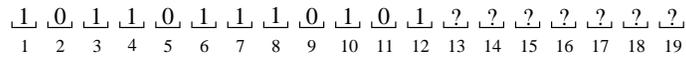, width=9cm}
        \caption{The first coordinates of $\text{s}(\text{s}(w))$}
        \label{fig:codage:2}
\end{center}
\end{figure}

As an example, consider Figure~\ref{fig:numerotation} (encoded by $w$) and
suppose an upward flip has been performed on 18 and 12. We know the values of
the normalized heights of $\text{s}(\text{s}(w))$ for vertices 1 to 12, as
shown in Figure~\ref{fig:codage:2}.

\begin{figure}[ht]
\begin{center}
        \epsfig{file=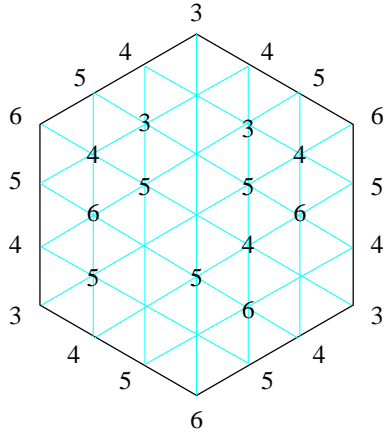,
                width=5cm}
        \caption{Utilization of the generalized Thurston algorithm}
        \label{fig:flip:TG}
\end{center}
\end{figure}

In order to use the generalized algorithm, we start with values as shown in
Figure~\ref{fig:flip:TG}.

\begin{figure}[ht]
\begin{center}
        \epsfig{file=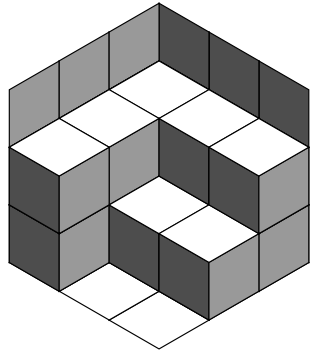,
                width=4cm}
        \caption{Result of the generalized Thurston algorithm}
        \label{fig:flip:resultat}
\end{center}
\end{figure}

\begin{figure}[ht]
\begin{center}
        \epsfig{file=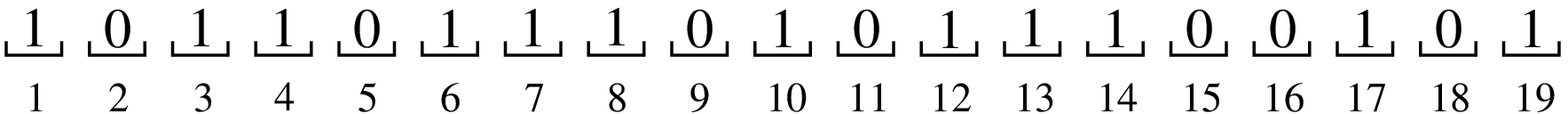, width=9cm}
        \caption{$\text{s}(\text{s}(w))$}
        \label{fig:codage:fin}
\end{center}
\end{figure}

The generalized algorithm yields the tiling of Figure~\ref{fig:flip:resultat}.
Once read, we obtain $\text{s}(\text{s}(w))$ (see
Figure~\ref{fig:codage:fin}).

We can now summarize the process in an algorithm:

\medskip

\noindent \textbf{\textit{Successor algorithm}}

\begin{itemize}

\item \textit{\textbf{Input}: A tileable domain $D$, a numerotation function
$\varphi$ from the set $V$ of the inner vertices of $D$ to $[\![ 1\,;\, |V|
|\!]$, the minimal tiling of $D$ and a word $w$ coding a tiling $T$ of $D$.}

\item \textit{\textbf{Step 1}: Compute the height function associated with
$T$ by using the normalized heights.}

\item \textit{\textbf{Step 2}: For $i$ from $|V|$ downto 1, examine whether
vertex $i$ is a local minimum of the height function; stop when such a vertex 
$i_0$ has been found. If no vertex is found, $w$ encodes the maximal tiling.}

\item \textit{\textbf{Step 3}: The first coordinates of $\text{s}(w)$ are:}

\[
	\begin{array}{l}
	\text{s}(w)[i] = w[i] \hskip1cm \text{for}\ 1 \leq i < i_0 \\
	\text{s}(w)[i_0] = w[i_0] + 1 \\
	\end{array}
\]

\item \textit{\textbf{Step 4}: Use the generalized Thurston algorithm in order
to find the smallest tiling bearing the heights of $\text{s}(w)$ for $i=1$ to
$i_0$\,.}

\item \textit{\textbf{Step 5}: Encode the tiling obtained at step~4 by a word,
which is $\text{s}(w)$.}

\item \textit{\textbf{Output}: $\text{s}(w)$.}

\end{itemize}

Let us analyze this algorithm. We denote by $|D|$ the number of vertices in
$D$.  Steps 1, 2, 3 and 5 require $O(|D|)$ operations; step~4
uses the generalized Thurston algorithm, which also runs in $O(|D|)$ time, so
the execution time is $O(|D|)$ and the algorithm is linear. The space required
is $O(|D|\:\ln|D|)$, which is the space needed to store a height function.

\soussection{Exhaustive generation}

Since the successor function preserves the lexicographic order, generating all
the tilings of $D$ amounts to recursively calling it:

\medskip

\noindent\textbf{\textit{Generation algorithm}}

\begin{itemize}

\item \textit{\textbf{Input}: A domain $D$.}

\item \textit{\textbf{Initialization}: Compute the minimal tiling of $D$ with 
Thurston's algorithm; stop if $D$ is not tileable. Encode the minimal tiling
by a word $w$. Numerote the vertices of $D$.}

\item \textit{\textbf{Current step}:} $w \longleftarrow \text{Successor}(w)$
\textit{and decode $w$ to obtain a tiling, until the} Successor
\textit{function does not produce a word.}

\item \textit{\textbf{Output}: the tilings of $D$.}

\end{itemize}

Let us analyze the algorithm. The time and space complexity are controlled by
the current step. Since the Successor function requires $O(|D|)$ operations,
the time complexity of the generation algorithm is $O(|D|)$ times the number
of tilings. In order to evaluate the space complexity, it is legitimate to
suppose that each tiling is discarded as soon as its successor is generated
since it is not used afterward. The space complexity is thus $O(|D|\:\ln|D|)$.

\section{Related algorithms}

The encoding of a tiling by a word, the successor function and the
generalization of Thurston's algorithm can be used to generate more than the
tilings of $D$: Indeed, all the characteristic elements of the lattice.

\soussection{Meet-irreducible elements}

The meet-irreducible elements of the lattice are those whose height function
admits exactly one local maximum in the interior of $D$. For each vertex in
the interior of $D$, one can compute the height in the minimal and maximal
tilings of $D$ using Thurston's original algorithm. The possible values for
$h(v)$ vary 3 by 3 in the case of lozenges, 4 by 4 in the case or dominoes, so
all the possibilities can easily be computed. For each pair defined by $v$ and
an admissible height, there exists a meet-irreducible element of the lattice,
which can be computed using the generalized Thurston algorithm.

\soussection{Order of the meet-irreducible elements}

The order on meet-irreducible elements is inherited from the lattice
structure. In the case of lozenges (resp. dominoes), such an element can
have at most 3 (resp. 2) successors in the order of the meet-irreducible
elements. Generating the full order thus amounts to examining whether the
putative successors are indeed meet-irreducible elements of the tiling.

\soussection{Lattice}

We have already generated all the tilings of $D$. In order to generate all the
lattice, it suffices to know which tilings can be obtained from a fixed tiling
$T$ by an upward flip.

This is done by examining the height function: for each vertex on which the
function admits a local minimum, there exists a single tiling that can be
deduced from $T$ by a single flip. Moreover, since the set of the tilings of
$D$ is connected by flips, all the links between tilings are thus obtained.

Another way to generate the lattice uses the order of meet-irreducible
elements: by Birkhoff's representation theorem, the former is isomorphic to
the order of the ideals of the latter and there exist optimal generic
algorithms that generate the order of the ideals of an arbitrary order
(see \cite{HMNS01} and \cite{KMNF92}).

\soussection{Intervals}

The lattice as a whole is a particular case of interval of itself. In order to
generate the elements of an interval, the successor function must be modified:
it is enough to compute the supremum of the result of the former successor
function and the minimal element of the interval. The links between the
tilings can, as above, be computed either by using the local minima of the
height function or by using the order of the meet-irreducible elements since
the intervals of a finite distributive lattice are themselves finite
distributive lattices.


\section{Conclusion}

We have provided several algorithms that make a non-trivial use of the lattice
structure of the tilings. Thurston's original algorithm has been
reinterpreted through Birkhoff's representation theorem and generalized in 
order to construct any tiling of the domain.

The normalized height functions provide a unified description of tilings by
dominoes and lozenges; they easily translate into a natural encoding of
tilings by words.

Our generation algorithm is linear in the number of tilings and requires a
space equivalent to the size of a single tiling. It can be extended to the
generation of all the characteristic elements of the lattice: meet-irreducible
elements and their order, the lattice structure and its intervals.

Furthermore, it should be rather straightforward to generalize the concepts
and methods to domains with holes.


\begin{thebibliography}{99}

\small

\bibitem[BH97]{Bur} J.K. Burton Jr, C.L. Henley, \textit{A constrained Potts
antiferromagnet model with an interface representation}, J. Phys. A {\bf 30}
(1997) 8385-8413.

\bibitem[Cha96]{Cha} T. Chaboud, \textit{Domino tiling in planar graphs with
regular and bipartite dual}, Theoretical Computer Science {\bf 159}
(1996), 137-142.

\bibitem[CL90]{CL90} J. H. Conway and J. C. Lagarias,
\textit{Tiling with Polyominoes and Combinatorial Group Theory},
Journal of Combinatorial Theory, A~53 (1990), p.~183-208.

\bibitem[Des01]{Des01} S. Desreux, \textit{An algorithm to
generate exactly once every tiling with lozenges of a domain}, to appear in
Theoretical Computer Science.

\bibitem[DP90]{Dav-Pri} P. A. Davey, H. A. Priestley, \textit{An introduction
to lattices and orders}, Cambridge University Press (1990).

\bibitem[HMNS01]{HMNS01} M. \textsc{Habib}, R. \textsc{Medina}, L.   
\textsc{Nourine}, G. \textsc{Steiner}, \emph{Efficient algorithms on   
distributive lattices}, Discrete Applied Mathematics~\textbf{110} (2001)   
169-187.

\bibitem[KMNF92]{KMNF92} T. \textsc{Kashiwabara}, S. \textsc{Masuda},
K. \textsc{Nakajima}, T. \textsc{Fujisawa}, \emph{Generation of maximum
independent sets of a bipartite graph and maximum cliques of a circular-arc
graph}, Journal of Algorithms~\textbf{13} (1992) 161-174.


\bibitem[Pro01]{Propp} J. G. Propp, \textit{Lattice Structure for Orientations
of Graphs}, preprint.

\bibitem[R\'em99]{Rem99} E. R\'emila, \textit{On the lattice structure of the
set of tilings of a simply connected figure with dominoes}, Proceedings of the
3rd International Conference on Orders, Algorithms and Applications (ORDAL)
(1999) LIP Research Rapport 1999-25, ENS Lyon.

\bibitem[Thu90]{Thu90} W. P. Thurston, \textit{Conway's Tiling Groups},
American Mathe\-ma\-ti\-cal Mon\-thly, \textbf{97}, oct.~1990, p.~757-773.






\end{thebibliography}
\end{document}